\documentclass[a4paper,11pt]{article}
\usepackage{amsmath,amsthm,amsfonts,amssymb,enumitem,graphicx,float,calc,microtype,mathtools,thmtools,underscore,anyfontsize,fontenc,lmodern}
\usepackage{latexsym}
\usepackage[dvipsnames]{xcolor}
\usepackage{todonotes}
 \usepackage{breakurl}
\usepackage{paralist}
\usepackage{enumitem}
\usepackage[english]{babel}
\usepackage[lmargin=30mm,rmargin=30mm,tmargin=30mm,bmargin=30mm]{geometry}
\renewcommand{\baselinestretch}{1.1}
\allowdisplaybreaks
\usepackage[utf8]{inputenc}
\usepackage[T1]{fontenc}

\newcommand{\defn}[1]{\textcolor{Maroon}{\emph{#1}}}
\usepackage[numbers,sort&compress,longnamesfirst]{natbib}
\def\NAT@spacechar{~}
\makeatother
\usepackage[pdftitle={Odd Colourings, Conflict-Free Colourings and Strong Colouring Numbers},
pdfauthor={Hickingbotham},
colorlinks =true,
linkcolor={black},
urlcolor={blue!80!black},
filecolor={blue!80!black},
citecolor={black},
anchorcolor=blue,
filecolor=blue, 
menucolor=blue]{hyperref}

\usepackage[noabbrev,capitalise]{cleveref}
\crefname{lem}{Lemma}{Lemmas}
\crefname{thm}{Theorem}{Theorems}
\crefname{cor}{Corollary}{Corollaries}
\crefname{prop}{Proposition}{Propositions}
\crefname{conj}{Conjecture}{Conjectures}
\crefname{open}{Open Problem}{Open Problems}
\crefname{obs}{Observation}{Observations}

\crefformat{equation}{(#2#1#3)}
\Crefformat{equation}{Equation #2(#1)#3}

\theoremstyle{plain}
\newtheorem{thm}{Theorem}

\theoremstyle{definition}

\renewcommand{\leq}{\leqslant}

\renewcommand{\geq}{\geqslant}

\theoremstyle{definition}

\DeclareMathOperator*{\bs}{\backslash}

\DeclareMathOperator*{\G}{\mathcal{G}}

\DeclareMathOperator{\scol}{scol}

\presetkeys%
{todonotes}%
{inline}{}

\begin{document}
	
	\title{\bf\Large Odd Colourings, Conflict-Free Colourings and Strong Colouring Numbers}
	\author{%
		Robert Hickingbotham\thanks{School of Mathematics, Monash University, Melbourne, Australia (\texttt{robert.hickingbotham@monash.edu}). Research supported by an Australian Government Research Training Program Scholarship.}
	}

	\date{\normalsize\today}
	\maketitle
	
	\begin{abstract}
		The odd chromatic number and the conflict-free chromatic number are new graph parameters introduced by Petru{\v{s}}evski and {\v{S}}krekovski [2021] and Fabrici, Lu{\v{z}}ar, Rindo{\v{s}}ov{\'a} and Sot{\'a}k [2022] respectively. In this note, we show that graphs with bounded $2$-strong colouring number have bounded odd chromatic number and bounded conflict-free chromatic number. This implies that graph classes with bounded expansion have bounded odd chromatic number and bounded conflict-free chromatic number. Moreover, it follows by known results that the odd chromatic number and the conflict-free chromatic number of $k$-planar graphs is $O(k)$ which improves a recent result of Dujmovi\'{c}, Morin and Odak [2022].
	\end{abstract}

	Let $G$ be a graph.\footnote{All graphs in this paper are finite, simple, and undirected. For $m,n \in \mathbb{Z}$ with $m \leq n$, let $[m,n]:=\{m,m+1,\dots,n\}$ and $[n]:=[1,n]$.}  A \defn{(vertex) $c$-colouring} of $G$ is any function $\psi: V(G) \to C$ where $|C|\leq c$. If $\psi(u)\neq \psi(v)$ for all $uv \in E(G)$, then $\psi$ is a \defn{proper colouring}. If $N(v):=\{w\in V(G):vw \in E(G)\}$ is the neighbourhood of a vertex $v$, then $\psi$ is an \defn{odd colouring} if for each $v \in V(G)$ with $|N(v)|>0$, there exists a colour $\alpha\in  C$ such that $|\{w \in N(v):\psi(w)=\alpha\}|$ is odd. Similarly, $\psi$ is a \defn{conflict-free colouring} of $G$ if for each $v \in V(G)$ with $|N(v)|>0$, there exists a colour $\alpha \in  C$ such that $|\{w \in N(v):\psi(w)=\alpha\}|=1$. The \defn{(proper) odd chromatic number} $\chi_{o}(G)$ of $G$ is the minimum integer $c$ such that $G$ has a (proper) odd $c$-colouring. Likewise, the \defn{(proper) conflict-free chromatic number} $\chi_{pcf}(G)$ of $G$ is the minimum integer $c$ such that $G$ has a (proper) conflict-free $c$-colouring. Clearly $\chi_{o}(G)\leq \chi_{pcf}(G)$ since a conflict-free colouring is an odd colouring.
	
	Motivated by connections to hypergraph colouring, the odd chromatic number and the conflict-free chromatic number were recently introduced by \citet{PS2021parity} and \citet{FLRS2022neighbourhoods} respectively. These parameters have gained significant traction with a particular focus on determining a tight upper bound for planar graphs. \citet{PS2021parity} showed that the odd chromatic number of planar graphs is at most $9$ and conjectured that their odd chromatic number is at most $5$. \citet{PP2022odd} improved this upper bound to $8$. For conflict-free colourings, \citet{FLRS2022neighbourhoods} proved a matching upper bound of $8$ for planar graphs. For proper minor-closed classes, a result of \citet{CLS2022odd1} implies that the odd chromatic number of $K_t$-minor free graphs is $O(t\sqrt{\log t})$. For non-minor closed classes, \citet{CLS2022odd1} showed that the odd chromatic number of $1$-planar graphs\footnote{A graph $G$ is \defn{$k$-planar} if it has an embedding in the plane such that each edge is involved in at most $k$ crossings.} is at most 23. \citet{DMO2022odd} proved a more general upper bound of $O(k^5)$ for the odd chromatic number of $k$-planar graphs. See \cite{CPS2022remarksodd,CPS2022remarksCF,CCKP2022oddsparse} for other results concerning these new graph parameters.
	
	In this note, we bound the conflict-free chromatic number by the $2$-strong colouring number. For a graph $G$, a total order $\preceq$ of $V(G)$, a vertex $v\in V(G)$, and an integer $s\geq 1$, let $R(G,\preceq,v,s)$ be the set of vertices $w\in V(G)$ for which there is a path $v=w_0,w_1,\dots,w_{s'}=w$ of length $s'\in[0,s]$ such that $w\preceq v$ and $v\prec w_i$ for all $i\in[s-1]$. For a graph $G$ and integer $s\geq 1$, the \defn{$s$-strong colouring number} $\scol_s(G)$ is the minimum integer $c$ such that there is a total order~$\preceq$ of $V(G)$ with $|R(G,\preceq,v,s)|\leq c$ for every vertex $v$ of $G$. 
	
	Colouring numbers provide upper bounds on several graph parameters of interest. First note that $\scol_1(G)$ equals the degeneracy of $G$ plus 1, implying $\chi(G)\leq \scol_1(G)$. A proper graph colouring is \defn{acyclic} if the union of any two colour classes induces a forest; that is, every cycle is assigned at least three colours. The \defn{acyclic chromatic number} $\chi_\text{a}(G)$ of a graph $G$ is the minimum integer $c$ such that $G$ has an acyclic $c$-colouring.
	\citet{KY2003orderings} proved that $\chi_\text{a}(G)\leq \scol_2(G)$ for every graph $G$. Other parameters that can be bounded by strong colouring numbers include weak colouring numbers \citep{zhu2009generalized}, game chromatic number \citep{KT1994uncooperative,KY2003orderings}, Ramsey numbers \citep{CS1993ramsey}, oriented chromatic number \citep{KSZ1997acyclic}, arrangeability~\citep{CS1993ramsey}, and boxicity \citep{EW2018boxicity}.
	
	Our key contribution is the following:
	\begin{thm}\label{MainResult}
		For every graph $G$, $\chi_{pcf}(G)\leq 2 \scol_{2}(G)-1$.
	\end{thm}
	
	Note that \cref{MainResult} is best possible as the conflict-free chromatic number is not bounded by the $1$-strong colouring number \cite{CPS2022remarksCF}. Before proving \cref{MainResult}, we highlight several noteworthy consequences. 
	
	First, \cref{MainResult} implies that graph classes with bounded expansion have bounded conflict-free colouring number and bounded odd chromatic number. Bounded expansion is a robust measure of sparsity with many characterisations \cite{zhu2009generalized,nevsetvril2008grad,nevsetvril2012sparsity}. For example, \citet{zhu2009generalized} showed that a graph class $\G$ has \defn{bounded expansion} if and only if there exists a function $f: \mathbb{N}\to \mathbb{N}$ such that $\scol_s(G)\leq f(s)$ for every graph $G \in \G$. Examples of graph classes with bounded expansion includes classes that have bounded maximum degree \cite{nevsetvril2008grad}, bounded stack number \cite{NOW2012examples}, bounded queue-number \cite{NOW2012examples}, bounded nonrepetitive chromatic number \cite{NOW2012examples}, or strongly sublinear separators \cite{DN2016sublinear}, as well as proper-minor closed classes \cite{nevsetvril2008grad}. See the book by \citet{nevsetvril2012sparsity} for further background on bounded expansion. As graph classes with bounded expansion have bounded $2$-strong colouring number, \cref{MainResult} implies that each of these classes have bounded conflict-free colouring number and thus have bounded odd chromatic number.
	
	Second, \cref{MainResult} implies a stronger bound for the odd chromatic number and the conflict-free chromatic number of $k$-planar graphs. Van den Heuvel and Wood \cite{HW2018improperARXIV} showed that $\scol_2(G)\leq 30(k+1)$ for every $k$-planar graph $G$. Thus we have the following consequence of \cref{MainResult}:
	\begin{thm}\label{kplanar}
		For every $k$-planar graph $G$, $\chi_{pcf}(G)\leq 60k+59.$
	\end{thm}
	
	\cref{kplanar} is the first known upper bound for the conflict-free chromatic number of $k$-planar graphs. For the odd chromatic number, the previous best known upper bound for $k$-planar graphs was $\chi_{o}(G)\in O(k^5)$ due to \citet{DMO2022odd}. 
	
	Finally, \cref{MainResult} gives the first known upper bound for the conflict-free chromatic number of $K_t$-minor free graphs. Van den Heuvel, Ossona de Mendez, Quiroz, Rabinovich and Siebertz \cite{HMQRS2017fixed} showed that $\scol_2(G)\leq 5\binom{t-1}{2}$ for every $K_t$-minor free graph $G$. Thus \cref{MainResult} implies the following:
	\begin{thm}\label{Ktminor}
		For every $K_t$-minor free graph $G$, $\chi_{pcf}(G)\leq  5(t-1)(t-2)-1.$
	\end{thm}
	 
	See \cite{HW2021shallow,HW2018improperARXIV,HMQRS2017fixed,DMN2021convex} for other graph classes that \cref{MainResult} applies to.

	\begin{proof}[Proof of \cref{MainResult}]
		We may assume that $G$ has no isolated vertices. Let $\preceq$ be the ordering $(v_1,\dots,v_n)$ of $V(G)$ where $|R(G,v_i,\preceq,2)| \leq \scol_{2}(G)$ for every vertex $v_i$ of $G$. For each vertex $v_i \in V(G)$, let $N^{-}(v_i):=R(G,v_i,\preceq,1)\setminus \{v_i\}$ be the \defn{left neighbours} of $v_i$, and let $v_j \in N(v_i)$ where $j=\min \{\ell \in [n]:v_{\ell}\in N(v_i)\}$ be the \defn{leftmost neighbour} of $v_i$. Let $\pi(v_i)$ denote the leftmost neighbour of $v_i$.
		
		We now specify the conflict-free colouring $\psi: V(G)\to [2 \scol_{2}(G)+1]$ by colouring the vertices left to right. For $i=1$, let $\psi(v_1)=1$. Now suppose $i>1$ and that $v_1,\dots,v_{i-1}$ are coloured. Let $X_i:=\{\psi(v_j):v_j\in R(G,u,\preceq,2)\setminus \{v_i\}\}$ and $Y_i:=\{\psi(\pi({v_j})):v_j\in N^{-}(v_i)\}$. Observe that $|X_i|\leq |R(G,u,\preceq,2)\setminus \{v_i\}|\leq \scol_{2}(G)-1$ and $|Y_i|\leq |R(G,v_i,\preceq,1)\setminus \{v_i\}|\leq \scol_{2}(G)-1$ and so $|X_i \cup Y_i|\leq 2\scol_{2}(G)-2$. As such, there exists some colour $\alpha\in [2 \scol_{2}(G)-1]\setminus (X_i \cup Y_i)$. Let $\psi(v_i):=\alpha$.
		
		Now $\psi$ is proper as each vertex receives a different colour to its left neighbours. We now show that it is conflict-free. Let $v_i\in V(G)$ and let $v_j=\pi(v_i)$. We claim that $\psi(v_j)\neq \psi(v_{\ell})$ for every $v_{\ell} \in N(v_i)\bs \{v_j\}$. Since $v_j$ is the leftmost neighbour of $v_i$, $j<\ell$. If $\ell < i$, then $v_j\in R(G,\preceq,v_{\ell},2)$ (by the path $v_{\ell},v_i,v_j$) and so $\psi(v_j)\in X_{\ell}$. Otherwise $i <\ell$ so $v_i \in N^{-}(v_{\ell})$ and thus $\psi(v_j)\in Y_{\ell}$. As such, $\psi(v_j)\in X_{\ell}\cup Y_{\ell}$ and hence $\psi(v_{j})\neq \psi(v_{\ell})$, as required.
	\end{proof}

	\subsection*{Acknowledgement}
	Thanks to David Wood for several helpful comments.

\fontsize{9.5}{10.5} 
\selectfont 
\let\oldthebibliography=\thebibliography
\let\endoldthebibliography=\endthebibliography
\renewenvironment{thebibliography}[1]{%
	\begin{oldthebibliography}{#1}%
		\setlength{\parskip}{0.2ex}%
		\setlength{\itemsep}{0.2ex}%
	}{\end{oldthebibliography}}
\bibliographystyle{DavidNatbibStyle}

\end{document}